\documentclass[a4paper,12pt,leqno]{article}
\usepackage{amsmath,amsfonts,amsthm,amssymb,dsfont}
\usepackage[alphabetic]{amsrefs}
\usepackage[OT4]{fontenc}

\long\def\symbolfootnote[#1]#2{\begingroup%
\def\thefootnote{\fnsymbol{footnote}}\footnote[#1]{#2}\endgroup}

\newtheorem{theorem}{Theorem}[section]
\newtheorem{proposition}[theorem]{Proposition}

\newtheorem{lemma}[theorem]{Lemma}

\theoremstyle{definition}
\newtheorem{construction}[theorem]{Construction}
\newtheorem{remark}[theorem]{Remark}
\newtheorem{definition}[theorem]{Definition}

\renewcommand{\proof}{\medskip\par\noindent\textbf{Proof.} \ignorespaces}
\renewcommand{\qed}{\quad\hskip0pt\null\hfill$\square$\par}

\def\Z{\mathbb{Z}}
\def\R{\mathbb{R}}

\def\C{\mathbb{C}}
\def\I{\mathrm{Id}}
\def\e{\vec{\mathbf{e}}}
\def\f{\vec{\mathbf{f}}}
\newcommand{\gl}{\mathrm{glue}}

\begin{document}

\begin{center}
\large\bfseries Veech groups of Loch Ness monsters
\end{center}

\begin{center}\bf
Piotr Przytycki$^a$\symbolfootnote[1]{Partially supported by MNiSW
grant N201 012 32/0718 and the Foundation for Polish Science.
}, Gabriela Schmith\"{u}sen$^b$\symbolfootnote[2]{Partially supported by Landesstiftung Baden-W\"{u}rttemberg.} \& Ferr\'an Valdez$^c$\symbolfootnote[3]{Partially supported by Sonderforschungsbereich/Transregio 45 and ANR Symplexe.}
\end{center}

\begin{center}\it
$^a$ Institute of Mathematics, Polish Academy of Sciences, \\
 \'Sniadeckich 8, 00-956 Warsaw, Poland \\
\emph{e-mail:}\texttt{pprzytyc@mimuw.edu.pl}
\end{center}

\begin{center}\it
$^b$
Institute of Algebra and Geometry, Faculty of Mathematics, \\
University of Karlsruhe, D-76128 Karlsruhe, Germany \\
\emph{e-mail:}\texttt{schmithuesen@mathematik.uni-karlsruhe.de}
\end{center}

\begin{center}\it
$^c$ Max-Planck Institut f\"{u}r Mathematik, \\
Vivatsgasse 7, 53111 Bonn, Germany \\
\emph{e-mail:}\texttt{ferran@mpim-bonn.mpg.de}
\end{center}

\begin{abstract}
\noindent
We classify Veech groups of tame non-compact flat surfaces. In particular we prove that all
countable subgroups of $\mathbf{GL}_+(2,\R)$ avoiding the set of mappings of norm less than $1$ appear as Veech groups of tame non-compact flat surfaces which are Loch Ness monsters. Conversely, a Veech group of any tame flat surface is either countable, or one of three specific types.
\end{abstract}



\section{Introduction}

For a compact flat surface $S$, the \emph{Veech group} of $S$ is the subgroup of $\rm \mathbf{SL}(2,\R)$ formed by the differentials of the orientation preserving affine homeomorphisms of $S$. Veech groups of compact flat surfaces are related to the dynamics of the geodesic flow  \cite{Ve}.

Our goal is to describe all possible Veech groups one can obtain for tame non-compact flat surfaces (see Definition \ref{flatsurface}), introduced in \cite{V2}. An example \emph{par excellence} of a tame non-compact flat surface is the surface associated to the billiard game on an irrational angled polygonal table. This surface is of infinite genus and has only one end \cite{V1}. A surface with those properties is called a \emph{Loch Ness monster} (see \cite{G}).
We distinguish the role of this "monster" in our main result.

\medskip
To state it, we need the following notation.
We denote by $\mathcal{U}\subset  \mathbf{GL}_+(2,\R)$ the set of matrices $M$ such that $||Mv||<||v||$ for all $v\in \R^2$, where $||\cdot||$ is the Euclidean norm on $\R^2$. We denote
\begin{itemize}
\item by $P\subset \mathbf{GL}_+(2,\R)$ the group of matrices
\begin{equation*}
\begin{pmatrix}
1 & t \\
 0 &s
\end{pmatrix}
, \mathrm{ where} \ t\in\mathbb{R}, s\in\mathbb{R}_+,
\end{equation*}
\item
by $P'\subset\mathbf{GL}_+(2,\R)$ the group of matrices generated by $P$ and $-\I$.
\end{itemize}
Note that $P$ has index $2$ in $P'$.

\medskip
We prove the following.

\begin{theorem}
\label{mtheorem1}
 Let $G\subset\mathbf{GL}_+(2,\R)$ be a Veech group of a tame flat surface. Then one of the following holds.
\begin{enumerate}
\item [(i)]$G$ is countable and disjoint from $\mathcal{U}$.
\item [(ii)] $G$ is conjugate to $P$.
\item [(iii)] $G$ is conjugate to $P'$.
\item[(iv)] $G=\mathbf{GL}_+(2,\mathbb{R})$.
\end{enumerate}
\end{theorem}

Conversely, we prove the following.

\begin{theorem}
\label{mtheorem2}
Any subgroup $G$ of $\mathbf{GL}_+(2,\R)$ satisfying assertion (i), (ii) or (iii) of Theorem \ref{mtheorem1} can be realized as a Veech group of a tame flat surface $X$ which is a Loch Ness monster.
\end{theorem}

In particular, every cyclic subgroup of $\mathbf{SL}(2,\R)$ or every Fuchsian group can be realized as the Veech group of a tame flat surface which is a Loch Ness monster. For compact flat surfaces, such questions are still open (see \cite[Problems 5, 6]{HZ}). Furthermore, observe that a cocompact Fuchsian group cannot be the Veech group of a compact flat surface \cite{Ve}, but occurs as the Veech group of a tame flat surface, which is a Loch Ness monster.

We will see that the only tame flat surfaces with Veech group $\mathbf{GL}_+(2,\R)$, as in (iv) of Theorem \ref{mtheorem1}, are cyclic branched coverings of the flat plane (see Lemmas \ref{one sing} and \ref{more sing}). In particular $\mathbf{GL}_+(2,\R)$ cannot be realized as a Veech group of a tame Loch Ness monster.

\medskip
In our article we restrict in Definition \ref{vgroup} of the Veech group to affine homeomorphisms which preserve the orientation. If we allow orientation reversing ones, substituting $\mathbf{GL}(2,\R)$ in place of $\mathbf{GL}_+(2,\R)$ in the statements of our theorems, they remain valid, except that we need to add three more ``parabolic'' subgroups to the pair $P$ and $P'$. No new ideas appear in the proofs. Thus we restrict to the orientation preserving case to simplify the formulation and the arguments.

\medskip
The article is organized as follows. In Section \ref{prel} we recall the definition
of a tame non-compact flat surface and its Veech group.

We divide the proofs of Theorems \ref{mtheorem1} and \ref{mtheorem2} into two parts.
In Section \ref{uncount} we treat the case where the group $G$ is uncountable. More precisely, we prove that if in the hypothesis of Theorem \ref{mtheorem1} we assume that $G$ is uncountable, then it satisfies assertion (ii), (iii) or (iv) (Proposition \ref{few possibilities}). Conversely, we prove that any group satisfying assertion (ii) or (iii) can be realized as a Veech group of a tame flat surface which is a Loch Ness monster (Lemmas \ref{example1} and \ref{example2}).

In Section \ref{count} we study the remaining case, where $G$ is countable. In other words, we prove that any group satisfying assertion (i) of Theorem \ref{mtheorem1} can be realized as a Veech group of a tame flat surface which is a Loch Ness monster (Proposition \ref{main constr}). This construction is the main point of the article. Conversely, we prove that if we assume in the hypothesis of Theorem \ref{mtheorem1} that $G$ is countable, then it satisfies assertion (i) (Lemma \ref{not in U}).

\medskip
\textbf{Acknowledgments}. We thank the faculty and staff of Max-Planck Institut in Bonn, where part of this work was carried out. We furthermore thank the Landesstiftung Baden--W\"urttemberg and the Department of Mathematics of the University of Karlsruhe that enabled the authors to meet and work together.

\section{Preliminaries}
\label{prel}

In this section we briefly recall the definition and features of non-compact flat surfaces. For more details, we refer the reader to \cite{V2}.

\medskip
Let $(S,\omega)$ be a pair formed by a connected Riemann surface $S$ and a non-zero holomorphic $1$--form $\omega$ on $S$. Denote by $Z(\omega)\subset S$ the zero locus of the form $\omega$. Local integration of $\omega$ endows $S\setminus Z(\omega)$ with an atlas whose transition functions are translations of $\C$. The pullback of the standard translation invariant flat metric on the complex plane defines a flat metric
on $S\setminus Z(\omega)$. Let $\widehat{S}$ be the metric completion of $S\setminus Z(\omega)$. Each point in $Z(\omega)$ has a neighborhood isometric to the neighborhood of $0\in \C$ with the metric coming from the 1--form $z^kdz$ for some $k>1$ (which is the metric induced via a cyclic branched covering of $\C$). The points in $Z(\omega)$ are called \emph{finite angle singularities}.

\begin{definition}
A point $p\in \widehat{S}$ is called an \emph{infinite angle singularity} of $S$, if there exists a neighborhood of $p$ isometric to the neighborhood of the branching point of the infinite cyclic branched covering of $\C$.
We denote the set of infinite angle
singularities of $\widehat{S}$ by $Y_{\infty}(\omega)$.
\end{definition}

\begin{definition}
\label{flatsurface}
The pair $(S,\omega)$ is called a \emph{tame flat surface}, if
$\widehat{S}\setminus S$ equals $Y_{\infty}(\omega)$.
\end{definition}

\medskip
Let $\mathrm{Aff}_+(S)$ be the group of affine orientation preserving homeomorphisms of a tame flat surface $S$ (we assume that $S$ comes with a preferred $1$--form $\omega$). Consider the differential
\begin{equation*}
\mathrm{Aff}_+(S)\overset{D}\longrightarrow \mathbf{GL}_+(2,\R)
\end{equation*}
that associates to every $\phi\in \mathrm{Aff}_+(S)$ its (constant) Jacobian derivative $D\phi$.

\begin{definition}
\label{vgroup}
Let $S$ be a tame flat surface. We call $G(S)=D(\mathrm{Aff}_+(S))$ the \emph{Veech group} of $S$.
\end{definition}

We define \emph{saddle connections} and \emph{holonomy vectors} in the context of tame non-compact flat surfaces exactly in the same way as for compact ones, see \cite{V2}.
\medskip

\medskip
We refer the reader to \cite{HS, Ve} for more details on Veech groups of compact flat surfaces, and to  \cite{HW, HSc, V2, H} for explicit examples of Veech groups of tame flat surfaces which are Loch Ness monsters.

\section{Uncountable Veech groups}
\label{uncount}

In this section we prove Theorems \ref{mtheorem1} and \ref{mtheorem2} in the case where $G$ is uncountable. Under this assumption we restate Theorem \ref{mtheorem1} in the following way.

\begin{proposition}
\label{few possibilities}
If the Veech group of a tame flat surface is uncountable, then it is conjugate to $P$, conjugate to $P'$ or equals the whole  $\mathbf{GL}_+(2,\R)$.
\end{proposition}

We begin the proof with the following.

\begin{lemma}
\label{one sing}
If a tame flat surface $S$ has no saddle connections and its Veech group $G$ is uncountable, then $G$ equals $P'$ or $\mathbf{GL}_+(2,\R)$. In the latter case $S$ is a cyclic branched covering of the flat plane.
\end{lemma}
\proof
First assume that $S$ has no singularities. Then the universal cover of $S$ is the flat plane and $S$ is either (i) the plane itself, or (ii) a flat cylinder which is a quotient of the plane by a cyclic group, or (iii) it is compact.
Since $G$ is uncountable, $S$ is not compact. In case (i) we have that $G=\mathbf{GL}_+(2,\R)$. In case (ii) we have that $G$ is conjugate to $P'$ by a rotation.

\medskip
Now assume that $S$ has a singularity $x_0$ (which might be of finite or infinite angle).
Since there are no saddle connections issuing from $x_0$, we have that $\widehat{S}$ is isometric to a (possibly infinite) cyclic branched covering of $\R^2$. Hence $G=\mathbf{GL}_+(2,\R)$.
\qed

\medskip
To complete the proof of Proposition \ref{few possibilities} it remains to prove the following.
\begin{lemma}
\label{more sing}
If the Veech group $G$ of a tame flat surface $S$ carrying saddle connections is uncountable, then $G$ is conjugate to $P$ or $P'$.
\end{lemma}
\medskip\par\noindent\textbf{Proof. Step 1.}\ignorespaces
\ \emph{All saddle connections of $S$ are parallel}.

\medskip
Since there are only countably many homotopy classes of arcs joining singularities of $\widehat{S}$, the set of saddle connections of $S$, and thus the set $V\subset \R^2$ of holonomy vectors, is countable. If $s_1$ and $s_2$ are two non-parallel saddle connections, then let $v_1$, $v_2$ be their holonomy vectors. For each $g\in G$ we define $\eta(g)=(g(v_1),g(v_2))\in V\times V$. Since $\{v_1,v_2\}$ is a basis of $\R^2$, we have that $\eta$ is an embedding. But $V\times V$ is countable. Contradiction. This concludes Step 1. 	

\medskip
Without loss of generality we may assume that all saddle connections are horizontal.
Let  ${\rm Spine}(S)\subset \widehat{S}$ be the union of the set of singularities together with all singular horizontal geodesics (this includes saddle connections). We claim that ${\rm Spine}(S)$ is connected and complete w.r.t. its intrinsic path metric. The latter follows from the completeness of $\widehat{S}$. The former follows from the fact that any two singularities of $\widehat{S}$ are connected by a concatenation of saddle connections, which are horizontal by Step 1.

\medskip\par\noindent\textbf{Step 2.}\ignorespaces
\ \emph{We have that $P\subset G$.}

\medskip
Let $C$ be the closure of a component of $\widehat{S}\setminus {\rm Spine}(S)$. It is a complete Riemann surface with nonvanishing holomorphic $1$--form and horizontal boundary. The boundary of $C$ is connected, since otherwise there would be a non-horizontal saddle connection joining singularities in different boundary components. Hence $C$ is either a half-plane or a half-cylinder with horizontal boundary. In particular, for any $g\in P$ we have that $C$ admits an orientation preserving affine homeomorphism with differential $g$, which fixes its boundary.
Hence for any $g\in P$, there is an orientation preserving affine homeomorphism $\overline{g}\in{\rm Aff_+}(S)$, with $D\overline{g}=g$, which fixes ${\rm Spine}(S)$ and is prescribed independently on each component of the complement.

\medskip\par\noindent\textbf{Step 3.}\ignorespaces
\ \emph{We have that $G\subset P'$.}

\medskip
Let $\e$ denote the unit horizontal vector in $\R^2$. We prove that for every $g\in G$ we have $g(\e)=\pm\e$. Otherwise, assume that there is an orientation preserving affine homeomorphism $\overline{g}\in{\rm Aff_+}(S)$ with differential $g$ for which $g(\e)=\lambda\e$, with $|\lambda|\neq 1$. Then $\overline{g}$ or its inverse acts as a contraction on ${\rm Sing}(S)$. By the Banach fixed point theorem, the iterates of any singularity under $\overline{g}$ or its inverse accumulate on the fixed point of $\overline{g}$. Since the set of singularities is invariant under the action of $\overline{g}$, this implies that it has an accumulation point. Contradiction.

\medskip
We summarize. By Steps 2 and 3 we have that $P\subset G\subset P'$. Since $P$ is of index $2$ in $P'$, we have that $G=P$ or $G=P'$.
\qed

\medskip
We now provide examples of Loch Ness monsters with Veech groups $P$ and $P'$. First we introduce the following vocabulary, which will become particularly useful in Section \ref{count}.

\begin{definition}
Let $S$ be a tame flat surface. A \emph{mark} on $S$ is an oriented finite length geodesic (with endpoints) on $S$ which does not meet singularities. If $S$ is simply connected, a mark is determined by its endpoints. The \emph{slope} of a mark is its holonomy vector, which lies in $\R^2$.

If $m, m'$ are two disjoint marks on $S$ with equal slopes, we can perform the following operation. We cut $S$ along $m$ and $m'$, which turns $S$ into a surface with boundary consisting of four straight segments.  Then we reglue these segments to obtain a tame flat surface $S'$ different from the one we started from. We say that $S'$ is obtained from $S$ by \emph{regluing along $m$ and $m'$}.

Let $S_0=S\setminus (m\cup m')$. Then $S'$ admits a natural embedding $i$ of $S_0$. If $A\subset S_0$, then we say that $i(A)$ is \emph{inherited} by $S'$ from $A$.
\end{definition}

\begin{remark}
\label{Euler and 4pi}
If $S'$ is obtained from $S$ by regluing, then
the number of singularities of $S'$ of a fixed angle equals the one of $S$, except for $4\pi$--angle singularities, whose number is greater by $2$ in $S'$ (we put $\infty+2=\infty$). The Euler characteristic of $S$ is greater by $2$ than the Euler characteristic of $S'$.
\end{remark}

We can extend the notion of regluing to families of marks.

\begin{definition}
Let $S$ be a tame flat surface. Assume that $\mathcal{M}=(m_n)_{n=1}^{\infty}$ and $\mathcal{M'}=(m'_n)_{n=1}^{\infty}$ are ordered families of disjoint marks, which do not accumulate in $\widehat{S}$, and such that the slope of $m_n$ equals the slope of $m'_n$, for each $n$. Let $S_0=S$ and let $S_n$ be obtained from $S_{n-1}$ by regluing along $m_n$ and $m'_{n}$. Let $S'$ be the Riemann surface equipped with a holomorphic $1$--form which is the limit of $S_n$. The limit exists since the marks do not accumulate, but might not be a tame flat surface. We say that $S'$ is obtained from $S$ by \emph{regluing along $\mathcal{M}$ and $\mathcal{M}'$}. If $A\subset S\setminus (\mathcal{M}\cup \mathcal{M'})$, then we define the subset of $S'$ \emph{inherited} from $A$ as before.
\end{definition}

We are ready to perform the following constructions.

\begin{lemma}
\label{example1}
There is a tame Loch Ness monster with Veech group $P$.
\end{lemma}
\proof
Let $A$ and $A'$ be two oriented flat planes, equipped with origins that allow us to identify them with $\R^2$. Let $\mathcal{C},\mathcal{C}'$ be families of marks with endpoints $(4n+1)\e, (4n+3)\e$, for $n\geq 1$, on $A,A'$, respectively, where $\e$ denotes, as before, the horizontal unit vector in $\R^2$. Let $\hat{A}$ be the tame flat surface obtained from $A\cup A'$ by regluing along $\mathcal{C}$ and $\mathcal{C'}$.

The group $P$ acts on $A$ and $A'$ under identification with $\R^2$. This action carries over to $\hat{A}$. Hence the Veech group $G$ of $\hat{A}$ contains $P$. By Lemma \ref{more sing}, we have that $G=P$ or $G=P'$. But in the latter case, the affine homeomorphism with differential $-\I$ must act on ${\rm Sing}(\hat{A})$ (defined in the proof of Lemma \ref{more sing}) by an orientation reversing isometry. Since there is no such isometry, we conclude that $G=P$.

By Remark \ref{Euler and 4pi}, we have that $\hat{A}$ has infinite genus. It has one end (this follows in particular from Lemma \ref{one end}). Hence $\hat{A}$ is a Loch Ness monster with Veech group $P$.
\qed

\begin{lemma}
\label{example2}
There is a tame Loch Ness monster with Veech group $P'$.
\end{lemma}
\proof
Similarly as in the proof of Lemma \ref{example1}, let $A$ and $A'$ be two oriented flat planes, equipped with origins that allow us to identify them with $\R^2$. Let $\mathcal{C},\mathcal{C}'$ be families of marks with endpoints $(4n+1)\e, (4n+3)\e$, on $A,A'$, respectively, where this time we take $n\in \Z$, and we order the marks into sequences. Let $\hat{A}$ be the tame flat surface obtained from $A\cup A'$ by regluing along $\mathcal{C}$ and $\mathcal{C'}$.

This time the action of the whole group $P'$ carries over to $\hat{A}$. Hence the Veech group $G$ of $\hat{A}$ contains $P'$. By Lemma \ref{more sing} we have that $G=P'$. The surface $\hat{A}$ is a Loch Ness monster by the same argument as in the proof of Lemma \ref{example1}.
\qed

\medskip
Lemmas \ref{example1} and \ref{example2} prove Theorem \ref{mtheorem2} in the case where $G$ is uncountable.

\section{Countable Veech groups}
\label{count}
The main part of this section is devoted to the proof of Theorem \ref{mtheorem2} in the case where the group $G\subset \mathbf{GL}_+(2,\R)$ is countable. In other words, we prove the following.

\begin{proposition}
\label{main constr}
For any countable subgroup $G$ of $\mathbf{GL}_+(2,\R)$ disjoint from $\mathcal{U}=\{g\in  \mathbf{GL}_+(2,\R)\colon ||g||<1\}$ there exists a tame flat surface $S=S(G)$, which is a Loch Ness monster, with Veech group $G$.
\end{proposition}
In fact the group $\mathrm{Aff}_+(S)$ will map isomorphically onto $G$ under the differential map. This means that the group $G$ will act on $S$ via affine homeomorphisms with appropriate differentials. Here we adopt the convention that an action of a group $G$ on a set $X$ is a mapping $(g,x)\rightarrow g\cdot x$ such that $(gh)\cdot x=g\cdot (h\cdot x)$ and $\mathrm{Id}\cdot x=x$.

\medskip
We begin with an outline of the proof of Proposition \ref{main constr}. We make use of the fact that any group $G$ acts on its Cayley graph $\Gamma$. We turn $\Gamma$ equivariantly into a flat surface. With each vertex $g$ of $\Gamma$ we associate a flat surface $V_g$ which can be cut into a flat plane $A_g$ and a \emph{decorated surface} $\widetilde{L}'_g$, whose role is explained later.

To guarantee tameness, we do not want the singularities of different $V_g$ to accumulate. Let $(g,g')$ be an edge of $\Gamma$ such that $g^{-1}g'$ is the $i$'th generator of $G$. We associate to this edge a \emph{buffer surface} $\hat{E}^i_g$ which connects $V_g$ to $V_{g'}$, but separates them by a definite distance.

We keep track of the end in the following way. First we provide that each $V_g$ and $\hat{E}^i_g$ is one-ended. Then we provide that after gluing all $V_g$ and $\hat{E}^i_g$, their ends actually merge into one end.

In this way we construct a one-ended flat surface with a faithful affine action of $G$. The role of the decorated surface $\widetilde{L}'_g$ is to prevent the group of orientation preserving affine homeomorphisms of the surface from being richer than $G$. To achieve this, $\widetilde{L}'_g$ is decorated with special singularities. This guarantees that every orientation preserving affine homeomorphism of the surface permutes this set of singularities and with some more care we establish that it actually acts as one of the elements of $G$.

\medskip
We begin by explaining how to obtain a nice action of $\mathbf{GL}_+(2,\R)$ on a disjoint union of affine copies of any flat surface.

\begin{definition}
\label{composing}
Let $S_{\I}$ be a tame flat surface. For each $g\in \mathbf{GL}_+(2,\R)$, we denote by $S_g$ the affine copy of $S_{\I}$, whose atlas differs from the one of $S_{\I}$ by post-composing each chart with $g$. In other words, $S_g$ comes with a canonical affine homeomorphism $\overline{g}\colon S_{\I}\rightarrow S_g$ with differential $g$.
Moreover, $\mathbf{GL}_+(2,\R)$ acts on the union of all $S_{g'}$ so that $\overline{g}$ maps each $S_{g'}$ onto $S_{gg'}$, with differential $g$.
\end{definition}

We provide the following criterion for $1$--endedness. Let $\Gamma$ be a connected graph. Let $A$ be the union, over $v\in \Gamma^{(0)}$, of $1$--ended tame flat surfaces $A_v$ without infinite angle singularities. Assume that each $A_v$ is equipped with infinite families of marks $\mathcal{C}^e_v$, for each edge $e$ issuing from $v$, and additional, possibly finite, two families of marks $\mathcal{C}_v, \mathcal{C}'_v$, of the same cardinality. Assume that all these marks are disjoint and do not accumulate. In particular this implies that $\Gamma^{(0)}$ is countable. Moreover, assume that for each edge $e=(v,v')$ the slopes of the marks in $\mathcal{C}^e_v$ and $\mathcal{C}^e_{v'}$ agree. Additionally, assume that the slopes of the marks in $\mathcal{C}_v$ and $\mathcal{C}'_v$ agree.

\begin{lemma}
\label{one end}
Let $S$ be the surface obtained from $A$ by regluing along $\mathcal{C}^e_v$ and $\mathcal{C}^e_{v'}$, for all edges $e=(v,v')$ in $\Gamma^{(1)}$, and along $\mathcal{C}_v$ and $\mathcal{C}'_v$, for all vertices $v$ in $\Gamma^{(0)}$. Then $S$ is $1$--ended. If $\Gamma$ has an edge or if it has only one vertex $v$ but with infinite $C_v$ (or if $A_v$ has infinite genus), then $S$ has infinite genus.

Unless $\Gamma$ has no edges (it has then only one vertex $v$) and additionally $\mathcal{C}_v$ is finite and $A_v$ has finite genus, we have that $S$ has infinite genus.
\end{lemma}
\proof
For each vertex $v$ in $\Gamma^{(0)}$, choose a basepoint $O_v$ in $A_v$. Let $B_v(r)$ be the closure in $S$ of the subset inherited from the ball of radius $r$ around $O_v$ with appropriate marks removed.

We order all vertices of $\Gamma$ into a sequence $(v_j)_{j=1}^{\infty}$. For $l\geq 1$, let $$K_{l}=\bigcup_{j=1}^{l} B_{v_j}(l).$$
Then $K_{l}$ is a family of compact sets which has the property that each compact set in $S$ is contained in $K_{l}$, for some $l\geq 1$.

Now we prove that the complement of each $K_{l}$ is connected. Since the $A_v$ are complete non-positively curved and $1$--ended, since balls and the marks we consider are convex, and since those marks are disjoint, we have that all $$A_{v_j}'=A_{v_j}\setminus (B_{v_j}(l)\cup_e \mathcal{C}_{v_j}^e\cup \mathcal{C}_{v_j}\cup \mathcal{C}'_{v_j})$$ are connected. Since $\Gamma$ is connected, all $\mathcal{C}^e_v$ are infinite, and $K_l$ intersects only a finite number of marks, we have that all $A'_{v_j}$ are in the same connected component of $S\setminus K_l$. Since the union of $A'_{v_j}$ is dense in $S\setminus K_l$, this implies that $S\setminus K_l$ is connected.

Thus $S$ is $1$--ended. If $\Gamma$ has at least one edge or $C_v$ is infinite, then $S$ has infinite genus by Remark \ref{Euler and 4pi}.
\qed

\medskip
We describe the construction of the \emph{buffer surfaces}, which will correspond to the edges of the Cayley graph $\Gamma$ of $G$. We denote the base vectors $(1,0),(0,1)$ of $\R^2$ by $\e$ and $\f$, respectively.

\begin{construction}
\label{buffer}
Let $E_{\I}, E'_{\I}$ be two oriented flat planes, equipped with origins that allow us to identify them with $\R^2$. We define the following families of slope $\e$ marks on $E_{\I}\cup E'_{\I}$. Let $\mathcal{S}$ be the family of marks on $E_{\I}$ with endpoints $4n\e, (4n+1)\e$, for $n\geq 1$, and let $\mathcal{S}_\gl$ be the family of marks on $E_{\I}$ with endpoints $(4n+2)\e, (4n+3)\e$, for $n\geq 1$. Let $\mathcal{S}'$ be the family of marks on $E'_{\I}$ with endpoints $2n\f, 2n\f+\e$, for $n\geq 1$, and let $\mathcal{S}'_\gl$ be the family of marks on $E'_{\I}$ with endpoints $(2n+1)\f,(2n+1)\f+\e$, for $n\geq 1$. Let $\hat{E}_{\I}$ be the tame flat surface obtained from $E_{\I}$ and $E'_{\I}$ by regluing along $\mathcal{S}_\gl$ and $\mathcal{S}'_\gl$. We call $\hat{E}_{\I}$ the \emph{buffer surface}. We record that $\hat{E}_{\I}$ comes with distinguished families of marks inherited from $\mathcal{S},\mathcal{S}'$, for which we retain the same notation.
\end{construction}

\begin{lemma}
\label{buffer lemma}
Let $\hat{E}_{\I}$ be the buffer surface and let $g\in \mathbf{GL}_+(2,\R)\setminus \mathcal{U}$. Then the distance in $\hat{E}_g$ (see Definition \ref{composing}) between $\overline{g}\mathcal{S}$ and $\overline{g}\mathcal{S}'$ is at least $\frac{1}{\sqrt{2}}$.
\end{lemma}

\proof
Denote by $\hat{d}$ the distance in $\hat{E}_g$ between $\overline{g}\mathcal{S}$ and $\overline{g}\mathcal{S}'$. Let $d$ be the distance in $E_g$ between $\overline{g}\mathcal{S}$ and $\overline{g}\mathcal{S}_\gl$ and let $d'$ be the distance in $E'_g$ between $\overline{g}\mathcal{S}'_\gl$ and $\overline{g}\mathcal{S}'$. Then we have that $\hat{d}\geq d+d'$.
Moreover, $d=|g(\e)|$ and $$d'=\min_{|s|\leq 1} |g(\f+s\e)|.$$ Let $s\in [-1,1]$ be such that the minimum is attained, that is $d'=|g(\f+s\e)|$. If $d+d'<\frac{1}{\sqrt{2}}$, then
$$|g(\f)|\leq |g(\f+s\e)|+|s||g(\e)|< \frac{1}{\sqrt{2}}.$$ Hence for any $v=x\e+y\f\in \R^2$ we have that
$$|g(v)|\leq |x||g(\e)|+|y||g(\f)|< \frac{1}{\sqrt{2}}(|x|+|y|)\leq \sqrt{x^2+y^2}=|v|.$$ Thus $||g||<1$. Contradiction.
\qed

\medskip
Now we construct the \emph{decorated surface} which will force rigidity of the affine homeomorphism group.

\begin{construction}
\label{decorated surface}
Let $L_{\I}$ be an oriented flat plane, equipped with an origin. Let $\widetilde{L}_{\I}$ be the threefold cyclic branched covering of $L_{\I}$, which is branched over the origin.
Denote the projection map from $\widetilde{L}_\I$ onto $L_\I$ by $\pi$. Denote by $R$ the closure in $\widetilde{L}_{\mathrm{Id}}$ of one connected component of the pre-image under $\pi$ of the open right half-plane in $L_{\mathrm{Id}}$. On $R$ consider coordinates induced from $L_{\mathrm{Id}}$ via $\pi$. Denote by $\mathcal{C}'$ the family of marks in $R$ with endpoints $(2n-1)\e, 2n\e$, for $n\geq 1$, and denote by $t$ and $b$ the two marks in $\widetilde{L}_{\mathrm{Id}}$ with endpoints in $R$ with coordinates $\f,2\f$ and $-\f, -2\f$, respectively. Let $\widetilde{L}'_{\I}$ be the tame flat surface obtained from $\widetilde{L}_{\I}$ by regluing along $t$ and $b$. We call $\widetilde{L}'_{\I}$ the \emph{decorated surface}.
\end{construction}

\begin{remark}
\label{position of O}
We keep the notation $\mathcal{C}'$ for the family of marks inherited by $\widetilde{L}'_{\I}$. We denote the point inherited from the origin by $O$. Then $O$ is a $6\pi$--angle singularity outside $\mathcal{C}'$.
\end{remark}

\begin{remark}
\label{3 saddle connections}
Let $S$ be a tame flat surface with a non-accumulating (in $\widehat{S}$) family $\mathcal{C}$ of marks with slopes $\e$. Assume that $S'$ is obtained from $\widetilde{L}'_\I\cup S$ by regluing along $\mathcal{C}'$ and $\mathcal{C}$. Then
there are only three saddle connections issuing from the point inherited from $O$ by $S'$. Their interiors are all contained in the subset inherited from $R\setminus (t\cup b\cup \mathcal{C}')$ and their holonomy vectors equal $-\f, \e$, and $\f$. Hence the angles between these saddle connections are $\frac{\pi}{2}, \frac{\pi}{2}$ and $5\pi$.
\end{remark}

\medskip
We are now ready for our main construction.
Recall that $\mathcal{U}$ denotes the set of linear mappings of norm less than one.
\begin{construction}
\label{monster}
Let $G$ be a nontrivial countable subgroup of $\mathbf{GL}_+(2,\R)\setminus\mathcal{U}$. Denote the generators of $G$ by $a_i$, where $i\geq 1$. If $G$ is trivial, we consider a single generator $a_1=\I$.
Let $A_{\I}$ be an oriented flat plane, equipped with an origin. Let $A$ be the union of $A_g$ over $g\in G$ (see Definition \ref{composing}).

\medskip
For $i\geq 0$ let $\mathcal{C}^i$ be the family of marks on $A_{\mathrm{Id}}$ with endpoints $i\f+(2n-1)\e,\ i\f+2n\e$, for $n\geq 1$. All these marks are pairwise disjoint.
Now, given $x_1, y_1\in \R$, consider the family $\mathcal{C}^{-1}$ of marks on $A_{\mathrm{Id}}$ with endpoints $(nx_1,y_1),\ (nx_1, y_1)+a_1^{-1}(\e)$, for $n\geq 1$. Choose $x_1>0$ sufficiently large and $y_1<0$ sufficiently small (i.e. $-y_1>0$ sufficiently large) so that all these
marks are pairwise disjoint and disjoint from the ones in $\mathcal{C}^i$ for $i\geq 0$.

Observe that a translate of the lower half-plane in $A_{\mathrm{Id}}$ is avoided by all already constructed marks.
In this way we can inductively, for all $i\geq 2$, choose $x_i,-y_i\in \R$ sufficiently large so that the marks with endpoints $(nx_i,y_i),\ (nx_i, y_i)+a_i^{-1}(\e)$, for $n\geq 1$, are pairwise disjoint and disjoint with the previously constructed marks. We denote these families by $\mathcal{C}^{-i}$. None of the described marks accumulate.

\medskip
Let $\widetilde{L}'_{\I}$ be the decorated surface from Construction \ref{decorated surface} and let $\widetilde{L}'$ be the union of $\widetilde{L}'_{g}$ over $g\in G$ (see Definition \ref{composing}).
For each $g\in G$ let $V_g$ be the flat surface obtained from $A_g\cup\widetilde{L}'_g$ by regluing along the families of marks $\overline{g}\mathcal{C}^0$ and $\overline{g}\mathcal{C}'$. The regluing is allowed, since all the slopes equal $g(\e)$. The surface $V_g$ is complete, in particular it is tame. Let $V$ be the union of the $V_g$ over $g\in G$. The action of $G$ on $A$ and on $\widetilde{L}'$ carries over to an action on $V$, and we retain the same notation for this action. It still has the property that the differential of $\overline{g}$ equals $g$, for each $g\in G$. We keep the notation $\mathcal{C}^i$, for $i\neq 0$, for the families of marks that are inherited from the families of marks on $A_{\I}$ by $V_{\I}$.

\medskip
For each $i\geq 1$ we consider a copy $\hat{E}^i_{Id}$ of the buffer surface $\hat{E}_{Id}$ defined in Construction \ref{buffer}. We denote the copies of $\mathcal{S},\mathcal{S}'$ in $\hat{E}^i_{Id}$ by $\mathcal{S}^i,\mathcal{S}'^i$. Let $E$ be the union of all $\hat{E}^i_g$, over $g\in G$ and all $i\geq 1$.
Let $S=S(G)$ be the Riemann surface equipped with the holomorphic $1$--form obtained from $V\cup E$ by regluing along the following pairs of families of marks. For each $i\geq 1$ and $g\in G$, we reglue the family $\overline{g}\mathcal{C}^i$ with $\overline{g}\mathcal{S}^i$ and the family $\overline{g}\mathcal{S}'^i$
with $\overline{g}\overline{a}_i\mathcal{C}^{-i}$. Note that this is allowed since all slopes of these marks equal $g(\e)$. Moreover, the action of $G$ carries over to $S$, and we retain the same notation for this action.
\end{construction}

\begin{remark}
\label{3pi singularities}
By Remarks \ref{Euler and 4pi} and \ref{position of O} the set of singularities of $S$ with angle $6\pi$ is the set of the $G$--translates of the point inherited by $S$ from $O$ (for which we retain the same notation). By Remark \ref{position of O} the translates $\overline{g}O$ of $O$ in $S$ are pairwise different, for different $g\in G$.
\end{remark}

\begin{lemma}
\label{lnm}
$S$ is a Loch Ness Monster.
\end{lemma}
\proof
This follows from Lemma \ref{one end} applied to the graph $\Gamma'$ obtained from the Cayley graph $\Gamma$ of $G=\langle a_i\rangle _{i\geq 1}$. We get $\Gamma'$ from $\Gamma$ by subdividing each edge of  $\Gamma$ into three parts and by adding for each original vertex $v$ of $\Gamma$ an additional vertex $v'$ and an edge joining $v'$ to $v$.
\qed

\begin{lemma}
\label{real surface}
$S$ is a tame flat surface.
\end{lemma}

\proof
Let $\bar{V}_g$, respectively $\bar{E}_g^i$, denote the closures in $S$ of the subsets inherited from
$V_g\setminus \overline{g}(\cup_{i\neq 0}\mathcal{C}^i)$, respectively $\hat{E}_g^i\setminus \overline{g}(\mathcal{S}^i\cup\mathcal{S}'^i)$.

It is enough to prove that $S$ is complete. Let $(x_k)$ be a Cauchy sequence on $S$. By Lemma \ref{buffer lemma} we may assume that there is some $g\in G$ such that all $x_k$ lie in the union of $\bar{V}_g$ and the adjacent affine buffer surfaces $\bar{E}^i_g$ and $\bar{E}^i_{ga^{-1}_i}$. Since the components of $\bar{V}_g\cap \left(\bigcup_i (\bar{E}^i_g\cup \bar{E}^i_{ga^{-1}_i})\right)$ form a discrete subset in $\bar{V}_g$, we may assume that all $x_k$ lie in $\bar{V}_g$ and in a single adjacent buffer surface. Since both $\bar{V}_g$ and the buffer surface are complete, $(x_k)$ converges, as required.
\qed

\begin{lemma}
\label{rigid}
Any orientation preserving affine homeomorphism of $S$ is equal to $\overline{g}$ for some $g\in G$.
\end{lemma}
\proof
Let $\psi$ be an orientation preserving affine homeomorphism of $S$. By Remark \ref{3pi singularities}, $\psi$ must permute the set of the $G$--translates of $O$. Hence $\psi(O)=\overline{g}(O)$, for some $g\in G$. We are going to prove that $\psi=\overline{g}$, which means that $\varphi=\overline{g}^{-1}\circ \psi$ equals the identity. For the time being we know only that $\varphi(O)=O$.

By Remark \ref{position of O}, there are only three saddle connections issuing from $O$. Exactly one angle formed by them at $O$ exceeds $\pi$. Hence $\varphi$, which is an orientation preserving affine homeomorphism fixing $O$, must fix all these saddle connections. Therefore $\varphi$ is equal to the identity in the neighborhood of $O$, which implies that $\varphi$ is the identity.
\qed

\medskip
We summarize with the following.

\medskip\par\noindent\textbf{Proof of Proposition \ref{main constr}.}\ignorespaces
\ If $G\subset \mathbf{GL}_+(2,\R)\setminus \mathcal{U}$ is countable, and nontrivial, then Construction \ref{monster} provides a Riemann surface $S=S(G)$ with a holomorphic $1$--form. Moreover, $G$ acts on $S$ by affine homeomorphisms with appropriate differentials. By Lemma \ref{real surface} the flat surface $S$ is tame. By Lemma \ref{lnm} it is a Loch Ness monster. By Lemma \ref{rigid} the Veech group of $S$ does not exceed $G$.
\qed

\medskip
This establishes Theorem \ref{mtheorem2} in the case where the group $G$ is countable.

\begin{remark}
\item[(i)]
If we do not require in Proposition \ref{main constr} that our flat surface is a Loch Ness monster, then it suffices to take only one mark from each infinite family of marks, instead of the whole family, in Construction \ref{monster}.
\item[(ii)]
If in Construction \ref{monster} we take, for positive odd $i$, the marks in $\mathcal{C}^i$ to have endpoints $i\f+(2n-1-\frac{1}{2^i})\e,\ i\f+(2n-\frac{1}{2^i})\e$, then there are Euclidean triangles of arbitrarily small area, with vertices in singularities, embedded in $S$. This is unlike in the case of compact flat surfaces, where small triangles appear only if the Veech group is not a lattice \cite{SW}.
\end{remark}

\medskip
Conversely, we have the following.

\begin{lemma}
\label{not in U}
If the Veech group $G$ of a flat surface $S$ is countable, then $G$ is disjoint from $\mathcal{U}$.
\end{lemma}
\proof
First consider the case, where $S$ has a singularity $x$.
Recall that $\widehat{S}$ denotes the metric completion of $S$ and that the action of the group of orientation preserving affine homeomorphisms of $S$ extends to an action on $\widehat{S}$. Suppose that there is an orientation preserving affine homeomorphism $\phi$ of $S$ with $D\phi\in \mathcal{U}$. Then $\phi$ extends to a contraction on $\widehat{S}$. By the Banach fixed point theorem, the sequence $\phi^k(x)$ converges in $\widehat{S}$. If $x$ is not the fixed point of $\phi$, then this contradicts tameness.

Assume now that $x$ is the fixed point of $\phi$ and the only singularity of $S$. Then $S$ is simply connected. Otherwise by pushing a homotopically nontrivial loop going through $x$ by the iterates of $\phi$ we obtain arbitrarily short homotopically nontrivial loops through $x$, which contradicts tameness. Hence $S$ is a cyclic branched covering of $\C$ and thus $G=\mathbf{GL}_+(2,\R)$ which is not countable, contradiction.

If $S$ does not have singularities, its universal cover is the flat plane. Since $G$ is countable, $S$ must be a flat torus and we have that $G\subset \mathbf{SL}(2,\R)$ which is disjoint from $\mathcal{U}$.
\qed

\medskip
This proves Theorem \ref{mtheorem1} in the case where $G$ is countable.


\begin{bibdiv}
\begin{biblist}


\bib{G}{article}{
   author={Ghys, {\'E}.},
   title={Topologie des feuilles g\'en\'eriques},
   language={French},
   journal={Ann. of Math. (2)},
   volume={141},
   date={1995},
   number={2},
   pages={387--422}
}

\bib{H}{article}{
  author={Hooper, P.},
  title={Dynamics on an infinite surface with the lattice property},
  date={2008},
  eprint={arXiv:0802.0189}
  }

\bib{HS}{article}{
   author={Hubert, P.},
   author={Schmidt, T. A.},
   title={An introduction to Veech surfaces},
   conference={
      title={Handbook of dynamical systems. Vol. 1B},
   },
   book={
      publisher={Elsevier B. V., Amsterdam},
   },
   date={2006},
   pages={501--526}
}

\bib{HSc}{article}{
author ={Hubert, P.},
author ={Schmith\"{u}sen, G.}
title  ={Infinite translation surfaces with infinitely generated Veech groups},
eprint ={http://www.cmi.univ-mrs.fr/~hubert/articles/hub-schmithuesen.pdf},
status ={preprint},
date={2008}
}

\bib{HW}{article}{
author ={Hubert, P.},
author ={Weiss, B.},
title  ={Dynamics on the infinite staircase surface},
eprint ={http://www.math.bgu.ac.il//~barakw/papers/staircase.pdf},
date={2008}
}

\bib{HZ}{article}{
   author={Hubert, P.},
   author={Masur, H.},
   author={Schmidt, T.},
   author={Zorich, A.},
   title={Problems on billiards, flat surfaces and translation surfaces},
   conference={
      title={Problems on mapping class groups and related topics},
   },
   book={
      series={Proc. Sympos. Pure Math.},
      volume={74},
      publisher={Amer. Math. Soc.},
      place={Providence, RI},
   },
   date={2006},
   pages={233--243}
}

\bib{SW}{article}{
author={Smillie, J.}
author={Weiss, B.}
title={Characterizations of lattice surfaces}
date={2008}
eprint={arXiv:0809.3729}
}

\bib{V1}{article}{
author ={Valdez, J. F.},
title  ={Infinite genus surfaces and irrational polygonal billiards},
  status={to appear},
  journal={Geom. Dedicata}
  date={2009}
}

\bib{V2}{article}{
author ={Valdez, J. F.},
title  ={Veech groups, irrational billiards and stable abelian differentials},
  eprint={arXiv:0905.1591v2},
  date={2009}
}

\bib{Ve}{article}{
   author={Veech, W. A.},
   title={Teichm\"uller curves in moduli space, Eisenstein series and an
   application to triangular billiards},
   journal={Invent. Math.},
   volume={97},
   date={1989},
   number={3},
   pages={553--583}
}

\end{biblist}
\end{bibdiv}
\end{document}